\documentclass[12pt]{article}
\usepackage{mathtools}
\usepackage{amsmath}
\usepackage[margin=1in]{geometry}
\usepackage{amsthm}
\usepackage{tikz}
\usepackage{amssymb}
\usepackage{tikz-qtree}
\usetikzlibrary{positioning}
\usepackage{lipsum}
\newtheorem{prop}{Proposition}[section]

\newtheorem{thm}[prop]{Theorem}

\theoremstyle{definition}
\newtheorem{eg}[prop]{Example}
\def \ds {\displaystyle} \def \tss {\textsuperscript} \def \p {\partial} \def \mc {\mathcal} \def \D {\Delta} \def \s {\sigma} \def \Var {\text{Var}} \def \mb {\mathbb} \def \a {\alpha} \def \r {\rho} \def \dz {\text{dz}} 
\begin{document}
\centerline
{
 \bf 
On Computing Certain Statistics on Ordered Trees
}

\bigskip
\centerline
{\it By Yonah BIERS-ARIEL}
\bigskip
{\bf Abstract}: We develop algorithms, implemented in a Maple package, that study the number of vertices with a particular number of children in a random ordered tree where all vertices must have a number of children in some finite set. By calculating the mixed moments of two such numbers, the package provides experimental evidence that the numbers are pairwise asymptotically normal.

\bigskip
{\bf Maple package and Sample Input and Output Files}: This article is accompanied by the Maple package {\tt ChildCountStatistics.txt} along with several input and output files available from the front of the article's webpage 
\\
{\tt http://sites.math.rutgers.edu/\~{}yb165/ChildCountStatistics/ChildCountStatistics.html}
\bigskip

\section{Overview}
We consider ordered, rooted trees where each vertex has a number of children in some finite set $S$, which we assume always contains 0 (since every tree has at least one leaf). For the remainder of the paper, anytime we refer to trees, we specifically mean this class of trees. In \cite{Andrew}, the authors consider the statistic $H_n$, the sum of the distances from each vertex to the root; here we perform similar analysis on the statistic $X_{n,s}$, the number of vertices with $s$ children in a tree with $n$ vertices overall.

In particular, we are interested in the moments of $X_{n,s}$. We would like to find $\mu_{n,s} = E[X_{n,s}]$ and $\sigma^2_{n,s}=\Var(X_{n,s})$, and, for higher powers $p$, we would like to find the scaled moment 
\[\frac{E[(X_{n,s}-\mu_{n,s})^p]}{(\s^2_{n,s})^{p/2}}.\]

Even better, we would like to take $X_{s_i}, X_{s_j}$ and find the $(p_i,p_j)$ scaled mixed moment 
 \[\frac{E[(X_{n,s_i}-\mu_{n,s_i})^{p_i}(X_{n,s_j}-\mu_{n,s_j})^{p_j}]}{(\s^2_{n,s_i})^{p_i/2} (\s^2_{n,s_j})^{p_j/2}}.\]
 
 From one perspective, this problem has already been solved. It is known that the distributions of child counts are asymptotically normal, and the scaled moments converge to the sequence $0,1,0,3,0,15,\dots$ of normal moments. Moverover, the distributions of child counts for different numbers of children are asymptotically jointly normal (see \cite{Janson} for example). From another perspective, however, virtually nothing is known about these moments since there is no easy means to compute them. The purpose of this paper and, more to the point, its accompanying Maple package, is to provide such a means.
  
 What we describe below is just one way to compute scaled moments; it is also possible to do so using the algebraic generating function ansatz. We believe our approach to be more efficient, though.

\section{Method}
The first step in this analysis is finding the total number of trees with $n$ vertices. Every tree with $n$ vertices is either a root with no children or else a root with some number of children, each of which is the root of a subtree. Algebraically, if $\mathcal{T}(S)$ is the set of all trees where each vertex has a number of children in $S$, we represent this as 
\begin{equation}\label{alg} \mathcal{T}(S) =\bigcup_{i \in S} \{\cdot\} \times \mathcal{T}(S)^i.\end{equation}

Let $f_n$ be the number of trees in $\mathcal{T}(S)$ with $n$ vertices, and define 
\[f(x) = \sum_{n=0}^\infty f_n x^n = \sum_{T \in \mathcal{T}(S)} x\tss{\# of vertices in $T$}.\]

Using Equation \ref{alg}, we obtain the following algebraic equation for $f$:
\[f(x) = x\Big(\sum_{i \in S} f(x)^i\Big).\]

The tricky part is extracting from this the coefficient of a particular $x^n$ (we will denote this coefficient as $[x^n]$). To find it, we use the Lagrange Inversion Theorem, specifically the following statement given in \cite{Zeilberger}:

\begin{thm}[Lagrange Inversion Theorem]
If $u(x)$ and $\Phi(z)$ are formal Laurent series satisfying $u(x) = x\Phi(u(x))$, then $[x^n]u(x) = (1/n) [z^{n-1}] \Phi(z)^n$.
\end{thm}

Taking $u(x) = f(x)$ and $\Phi(z) = \sum_{i \in S} z^i$, our problem of computing the coefficients of $f$ is reduced to the problem of computing the coefficients of powers of $\Phi$. This can be done very efficiently using the amazing Almkvist-Zeilberger algorithm described in \cite{Almkvist}, which gives a recurrence satisfied by $f_n$. With this recurrence in hand, it is simple to find $f_n$ for as large an $n$ as is desired. The following example is a concrete demonstration of the use of this technique.

\begin{eg}\label{main}Suppose $S = \{0,1,2\}$; then $f$ must satisfy $f(x) = x(1+f(x) +f(x)^2)$. If we want to know the number of trees with $n$ vertices, then we need to find $f_n=[x^{n}]f(x)$, which, by the Lagrange Inversion Theorem, is equivalent to finding $1/n[z^{n-1}](1+z+z^2)^{n}$. This, in turn, is equivalent to finding the residual of $1/n (1+z+z^2)^{n}/z^{n}$; and, using the Almkvist-Zeilberger algorithm, we find that this residual satisfies the recurrence 
\[f_n = \frac{3n(n-1)f_{n-2} +(n+2)(2n+3)f_{n-1}}{(n+1)(n-1)}.\]
Since $f_1 = f_2 = 1$, we can now use this recurrence to find any desired value of $f_n$.
\end{eg}

Unfortunately, this technique so far only counts the number of trees; since we are also interested in the number of vertices with a particular number of children, we need to add some modifications. Letting $S = \{s_1,s_2,...s_k\}$ (note that $s_1=0$), we generalize $f(x)$ to be 
\[f(x;y_{s_1},y_{s_2},...,y_{s_k}) = \sum_{T \in \mathcal{T}(S)} x\tss{\# of vertices in $T$}\prod_{i=1}^k y_{s_i}\tss{\# of vertices with $s_i$ children}.\]

Let $\D_s$ be the operator $y_s\frac{\p}{\p y_s}$, and consider the result of applying $\D_s$ to $f$. Each monomial in $f$ corresponds to some tree in $\mc{T}(S)$, and applying $\D_s$ to $f$ simply multiplies each of these monomials by the number of vertices with $s$ children in the corresponding tree. Therefore, $\ds\D_s f(x;y_{s_1},y_{s_2},...,y_{s_k}) \Big|_{y_{s_1}=y_{s_2}=...=y_{s_k}=1}$ is the generating function for the total number of vertices with $s$ children. Similarly, $\ds\D_s^p f(x;y_{s_1},y_{s_2},...,y_{s_k}) \Big|_{y_{s_1}=y_{s_2}=...=y_{s_k}=1}$ is the generating function for the sum of the $p\tss{th}$ power of the number of vertices over all trees with $s$ children.

This method will tell us as much as we could want to know about the limiting marginal distributions of individual $X_{n,s}$; in particular it gives an empirical proof that these random variables are asymptotically normal (although not necessarily independent). However, we are more interested in whether $X_{n,s_i}, X_{n,s_j}$ are asymptotically jointly normal, and to do that we will need to compute mixed moments as well. Fortunately, doing so is just a matter of applying $\D_{s_i}$ and $\D_{s_j}$ as demonstrated in Example \ref{cont}.

We will not be computing these mixed moments directly, though. Instead, we will first calculate the numerators of these moments, which we now define. 

For random variables $X_1, X_2$ which assign values to elements of a set $A$, define $N_{p}(X) = \sum_{a \in A} X(a)^p$ and $N_{p_1,p_2}(X_1,X_2) = \sum_{a \in A} X_1(a)^{p_1}X_2(a)^{p_2}$. In our case, $A$ is the set of trees with $n$ vertices and child counts in $S$, while the $X$s will be $X_{s_i}$ and $X_{s_j}$. Therefore, $N_0(X_{n,s_i}) = N_0(X_{n,s_j}) = N_{0,0}(X_{n,s_i},X_{n,s_j})$ is the number of trees with child counts in $S$ on $n$ vertices. Meanwhile, $N_{1,0}(X_{n,s_i},X_{n,s_j}) = N_{1}(X_{n,s_i})$ is the sum over all such trees of the number of vertices with $s_i$ children, and $N_{1}(X_{n,s_i})/N_{0} (X_{n,s_i})= \mb{E}[X_{n,s_i}]$. In general, 
\[\frac{N_{p_1,p_2}(X_{n,s_i},X_{n,s_j})}{N_{0,0}(X_{n,s_i},X_{n,s_j})} = \mb{E}[X_{s_i}^{p_1} X_{s_j}^{p_2}].\]

We are now ready to return to Example \ref{main} and compute many mixed moments.

\begin{eg}\label{cont}
We continue the analysis begun in Example \ref{main} by considering the set $S =\{0,1,2\}$. Recall that in the earlier example we found that the number of trees on $n$ vertices with all child counts in $S$ followed the recurrence 
\[f_n = \frac{3n(n-1)f_{n-2} +(n+2)(2n+3)f_{n-1}}{(n+1)(n-1)}.\]
with initial conditions $f_1=f_2=1$. Therefore, there are 593742784829 trees on 30 vertices with child counts in $S$. 

Now, let's find $\mb{E}[X_{30,0}] =\mu_0$. As promised, we will use the generating function $f(x;y_0,...,y_k)$, where $f_n$ is still the coefficient on $x^n$. Denote the coefficient on $x^n$ in $\D_s f$ by $f_{n,y_s}$, and let $a_{n,y_s}=f_{n,y_s}|_{y_0=y_1=...=y_k=1}$. Note that we can find the total number of leaves on all trees with 30 vertices by finding $N_{1}(X_{30,0})$.

Note that the new $f$ satisfies the functional equation $f(x;y_0,y_1,y_2) = x(y_0+y_1f(x;y_0,y_1,y_2)+y_2f(x;y_0,y_1,y_2)^2)$. The Lagrange inversion theorem tells us that to find $f_n$ we need the coefficient of $z^{n-1}$ in $\frac{(y_0+y_1 z+y_2z^2)^n}{n}$, and so we have
\begin{equation*}
\resizebox{\textwidth}{!}{$f_{n,y_0} = y_0\frac{\p}{\p y_0} \oint_{|z|=1} \frac{(y_0+y_1z+y_2z^2)^n}{nz^n} \dz= \oint_{|z|=1} y_0\frac{\p}{\p y_0} \frac{(y_0+y_1z+y_2z^2)^n}{nz^n} \dz=\oint_{|z|=1} \frac{(z^2y_2+zy_1+y_0)^ny_0}{nz^n} \dz$ .}
\end{equation*}
This second integral is in exactly the form we'd like in order to use the Almkvist-Zeilberger algorithm. Using it, and plugging in 1 for all $y_i$, we obtain the recurrence
\[N_1(X_{n,0}) = \frac{(2(n-2)+1)N_1(X_{n-1,0}) + 3(n-2)N_1(X_{n-2,0})}{n+1}\]
with $N_1(X_{1,0}) = N_1(X_{2,0}) = 1$. With this recurrence, we find that $N_1(X_{30,0})= 6186675630819$, and so $\mb{E}[X_0] = N_1(X_{30,0})/N_0(X_{30,0}) = 10.42$. In other words, just over 1/3 of the vertices on 30-vertex trees with child counts in $S$ is a leaf.

We can do the same calculations for $X_1$ to find that $N_1(X_{n,1})$ satisfies the recurrence
\[N_1(X_{n,1}) = \frac{(n-1)(2n-3)N_1(X_{n-1,1}) + 3(n-1)(n-2)N_1(X_{n-2,1})}{n(n-2)}\]
with initial conditions $N_1(X_{1,1})=0, N_1(X_{2,1})=1, N_1(X_{3,1})= 2$. We find that $N_1(X_{30,1}) = 6032675068061$, and so $\mb{E}[X_0] = N_1(X_{30,1})/N_0(X_{30,1}) = 10.16$. 

Finally suppose that we want to find $\mb{E}[X_{30,0}^2X_{30,1}^3] = N_{2,3}(X_{30,0},X_{30,1})/ N_{0,0}(X_{30,0},X_{30,1})$. For this we calculate:
\[ (y_0\frac{\p}{\p y_0})^2(y_1\frac{\p}{\p y_1})^3 \oint_{|z|=1} \frac{(y_0+y_1z+y_2z^2)^n}{nz^n} = \oint_{|z|=1} (y_0\frac{\p}{\p y_0})^2(y_1\frac{\p}{\p y_1})^3 \frac{(y_0+y_1z+y_2z^2)^n}{nz^n}.\]
Running Almkvist-Zeilberger again, we obtain $ N_{2,3}(X_{30,0},X_{30,1}) = 68622906286794431$, and so $\mb{E}[X_{30,0}^2X_{30,1}^3] = 115576.83$.
\end{eg}

\section{From Raw Moments to Scaled Moments}\label{scale}

Up until now, we have only been concerned with calculating $N_{i,j}$, but we now want to turn them into scaled moments. Recall that we can calculate the $p_1\tss{th}$ moment about the mean, denoted $m_{p_1}$, as follows:
\begin{align*} m_{p_1} &= \mb{E}[(X_1-\mu_1)^{p_1}] = \mb{E}\Big[\sum_{r=0}^{p_1} \binom{p_1}{r}(-1)^r \mu_1^rX_1^{p_1-r}\Big] = \sum_{r=0}^{p_1} (-1)^r \binom{p_1}{r} \mu_1^r \mb{E}[X_1^{p_1-r}]
\\
&=\sum_{r=0}^{p_1} (-1)^r \binom{p_1}{r} \Big(\frac{N_{1,0}}{N_{0,0}}\Big)^r \frac{N_{p_1-r,0}}{N_{0,0}}=\frac{1}{N_{0,0}^{p_1}}\sum_{r=0}^{p_1} (-1)^r \binom{p_1}{r}N_{1,0}^rN_{0,0}^{p_1-r-1}N_{p_1-r}.
\end{align*}
For the problem here, we need to extend this to find $m_{p_1,p_2}$ by calculating:
\begin{align*}
m_{p_1,p_2} &= \mb{E}[(X_1-\mu_1)^{p_1}(X_2-\mu_2)^{p_2}] 
\\
&= \mb{E}\Big[\Big(\sum_{r=0}^{p_1} \binom{p_1}{r}(-1)^r \mu_1^rX_1^{p_1-r}\Big)\Big(\sum_{t=0}^{p_2} \binom{p_2}{t}(-1)^t\mu_2^tX_2^{p_2-t}\Big)\Big]
\\
&=\sum_{r=0}^{p_1}\sum_{t=0}^{p_2}\binom{p_1}{r}(-1)^r \mu_1^r\binom{p_2}{t}(-1)^t\mu_2^t \mb{E}[X_1^{p_1-r}X_2^{p_2-t}]
\\
&=\sum_{r=0}^{p_1}\sum_{t=0}^{p_2}\binom{p_1}{r}(-1)^r \Big(\frac{N_{1,0}}{N_{0,0}}\Big)^r\binom{p_2}{t}(-1)^t\Big(\frac{N_{0,1}}{N_{0,0}}\Big)^t \frac{N_{p_1-r,p_2-t}}{N_{0,0}}
\\
&=\frac{1}{N_{0,0}^{p_1+p_2}}\sum_{r=0}^{p_1} \sum_{t=0}^{p_2} \binom{p_1}{r} \binom{p_2}{t} (-1)^{r+t} N_{1,0}^r N_{0,1}^t N_{0,0}^{p_1+p_2-r-t-1}N_{p_1-r,p_2-t}.
\end{align*}

Letting $X_1 = X_{s_1}$ and $X_2 = X_{s_2}$, we can use the techniques from the previous section to find $N_{i,j}$ for any $i,j$ we like, and so combining them with this formula allows us to calculate $m_{p_1,p_2}$. Once we do the calculation we find the scaled mixed moment
\[\a_{p_1,p_2} = \frac{m_{p_1,p_2}}{m_{2,0}^{p_1/2} m_{0,2}^{p_2/2}}.\]

To find the corresponding normal moment, we begin by calculating the correlation $\rho = \a_{1,1}$. Now, if two variables $X,Y$ are jointly normally distributed with correlation $\r$, their pdf is given by
\[f(x,y) = \frac{e^{-\frac{x^2}{2}-\frac{y^2}{2}+\r xy}\sqrt{1-\r^2}}{2\pi},\]
and it is straightforward to compute the $(p_1, p_2)$ moment from this function. We then compare it to the actual $(p_1,p_2)$ moment of $X_{s_1},X_{s_2}$ to determine how close to normal that distribution is.

\section{Maple Implementation}
This paper is accompanied by the Maple package {\tt ChildCountStatistics.txt} It is capable of performing all the calculations described in this paper, including: finding recurrences for raw moments, finding the values of raw moments, computing scaled mixed moments, and computing the scaled mixed moments of a bivariate normal distribution whose correlation is the same as that between the number of vertices with $s_1$ children and the number with $s_2$ children.

The interested reader should download ChildCountStatistics.txt from this paper's website, and then read it in Maple. The function Help() will list all of the functions provided in the package and Help(function1) will explain what function1 does and give an example of how to use it.

Readers who only want to see the sorts of results that the package can produce can instead look at the sample input and output files, also available on this paper's website.

\section{Acknowledgements}
The author would like to thank his advisor, Doron Zeilberger, for suggesting this topic to him as well as for directing him to many relevant resources and helping with the Maple implementation. Thanks as well to an anonymous referee for suggesting references regarding existing proofs of asymptotic normality.

\end{document}